\documentclass[12pt]{article}

\usepackage{mathrsfs}
\usepackage{amsmath}
\usepackage{amsfonts}
\usepackage{amssymb}

\usepackage{amsmath, amsthm, amsfonts, amssymb}
\def\endpf{\hbox{\vrule height1.5ex width.5em}}

\newcommand{\CC}{{\mathbb C}}

\newcommand{\HH}{{\mathbb H}}
\newcommand{\BB}{{\mathbb B}}

\newcommand{\ra}{\rightarrow}
\def \-{\bar}
\newcommand{\p}{\partial}
\newcommand{\w}{\tilde}

\newtheorem{theorem}{Theorem}[section]
\newtheorem{lemma}[theorem]{Lemma}
\newtheorem{corollary}[theorem]{Corollary}
\newtheorem{proposition}[theorem]{Proposition}

\newtheorem{remark}[theorem]{Remark}
\newtheorem{claim}[theorem]{Claim}
\date{}

\begin{document}

\title{\bf Rigidity for local holomorphic isometric embeddings from
${\BB}^n$ into ${\BB}^{N_1}\times\cdots \times{\BB}^{N_m}$ up to conformal factors}

\author{
Yuan Yuan \ \ and \ \ Yuan Zhang}

\vspace{3cm} \maketitle

\begin{abstract}
{\small In this article, we study local holomorphic isometric
embeddings from ${\BB}^n$ into ${\BB}^{N_1}\times\cdots
\times{\BB}^{N_m}$ with respect to the normalized Bergman metrics up to conformal factors.
Assume that each conformal factor is smooth Nash algebraic. Then
each component of the map is a multi-valued holomorphic map between
complex Euclidean spaces by the algebraic extension theorem derived
along the lines of Mok and Mok-Ng. Applying holomorphic continuation
and analyzing real analytic subvarieties carefully, we show that
each component is either a constant map or a proper holomorphic map
between balls. Applying a linearity criterion of Huang, we conclude
the total geodesy of non-constant components.}
\end{abstract}

\bigskip
\section{Introduction}
Write $\BB^n:=\{z\in {\CC}^n: |z|<1\}$ for the unit ball in $\CC^n$.
Denote by $ds^2_n$ the normalized Bergman metric on $\BB^n$ defined
as follows:
\begin{align}\label{metric}
&ds^2_{n} =\sum_{j,k\le
n}\frac{1}{(1-|z|^2)^2}\big{(}(1-|z|^2)\delta_{jk}+\-z_jz_k\big)dz_j\otimes
d\-z_k.
\end{align}
Let $U\subset \BB^n$ be a connected open subset.
Consider a holomorphic isometric embedding
\begin{equation}
F=(F_1,\ldots, F_m): U
\rightarrow \BB^{N_1}\times \cdots
\times \BB^{N_m}
\end{equation}
up to conformal factors
$\{\lambda(z,\-{z});\lambda_1(z,\-{z}),\cdots, \lambda_m(z,\-{z})\}$
in the sense that
 $$\lambda(z,\-{z}) ds^2_n=\sum_{j=1}^{m}\lambda_j(z,\overline{z})
F_j^*(ds^2_{N_j}).$$ Here and in what follows,
the conformal factors $\lambda(z,\-{z}), \ \lambda_j(z,\-{z})$
($j=1,\cdots,m$) are assumed to be positive 
smooth
Nash algebraic functions over $\CC^n$. One can in fact assume that $\lambda(z, \bar z)=1$, and replace $\lambda_j(z, \bar z)$ by $\frac{\lambda_j(z, \bar z)}{\lambda(z, \bar z)}$. Under such notation, $\lambda_j(z, \bar z)$ is assumed to be positive, smooth and Nash algebraic.
Moreover, for each $j$ with
$1\le j\le m$, $ds_{N_j}^2$ denotes the corresponding normalized
Bergman metric of $\BB^{N_j}$ and $F_j$ is a holomorphic map from
$U$ to $\BB^{N_j}$. We write $F_j=(f_{j,1},\ldots, f_{j,l},\ldots,
f_{j, N_j})$, where $f_{j,l}$ is the $l$-th component of $F_j$. In
this paper, we prove the following rigidity theorem:

\begin{theorem} \label{main1}
Suppose $n\ge 2$. Under the above notation and assumption, we then
have, for each $j$ with $1\le j\le m$, that either $F_j$ is a
constant map or $F_j$ extends to a totally geodesic holomorphic
embedding from $({\BB}^n,ds^2_n)$ into $({\BB}^{N_j},ds^2_{N_j})$.
Moreover, we have the following identity
$$\sum_{F_j\hbox{ is not a constant}}\lambda_j(z,\-{z})=\lambda(z,\-{z}).$$
\end{theorem}

In particular, when $\lambda_j(z, \bar z),\lambda(z, \bar z)$ are
positive constant functions, we have  the following rigidity result for
local isometric embeddings:

\begin{corollary} \label{main2}
Let
\begin{equation}
F=(F_1,\ldots, F_m): U\subset {\BB^n}
\rightarrow
\BB^{N_1}\times \cdots \times \BB^{N_m}
\end{equation}
be a local holomorphic isometric embedding in the sense that
 $$\lambda ds^2_n=\sum_{j=1}^{m}\lambda_j
F_j^*(ds^2_{N_j}).$$ Assume that $n\ge 2$ and
 $\lambda,\ \lambda_j$ are positive constants.
We then have, for each $j$ with $1\le j\le m$, that either $F_j$ is
a constant map or $F_j$ extends to a totally geodesic holomorphic
embedding from $({\BB}^n,ds^2_n)$ into $({\BB}^{N_j},ds^2_{N_j})$.
Moreover, we have the following identity
$$\sum_{F_j\hbox{ is not a constant}}\lambda_j=\lambda.$$
 \end{corollary}

Recall that a  function $h(z,\-{z})$ is called a Nash algebraic
function over $\CC^n$ if
 there is an irreducible polynomial $P(z,\xi,X)$ in $(z,\xi,X)\in
{\CC^n}\times  {\CC^n}\times {\CC}$ with $P(z,\-{z},h(z,\-{z}))\equiv 0$ over $\CC^n$.
We mention that a holomorphic map from ${\BB}^n$ into ${\BB}^N$ is a
 totally geodesic embedding with respect to the normalized Bergman
 metric if and only if there are a (holomorphic) automorphism $\sigma
 \in Aut({\BB}^n)$ and an automorphism $\tau \in Aut({\BB}^N)$ such
 that $\tau\circ F\circ \sigma(z)\equiv(z,0)$. Also, we mention that
 by the work of Mok \cite{Mo1}, the result in Corollary 1.2 does not
 hold anymore when $n=1$. (See also many examples and related
 classification results in the work of Ng (\cite{Ng1}).

 The study of the global extension and  rigidity problem for local
 isometric embedding was first carried out in a  paper of Calabi
 \cite{Ca}. After [Ca],  there appeared quite a few papers along
 these lines of research (see \cite{Um}, for instance). In 2003,
 motivated by problems from Arithmetic Algebraic Geometry, Clozel and
 Ullmo \cite{CU} took up again the problem by considering the
 rigidity problem for a local isometric embedding with a certain
symmetry from ${\BB}^1$ into ${\BB}^1\times \cdots \times{\BB }^1$.
More recently, Mok carried out a systematic study of this problem in
a very general setting. Many far reaching deep results have been
obtained by Mok and later by Ng and Mok-Ng. (See \cite{Mo1}
\cite{Mo2} \cite{MN}  [Ng1-3] and the references therein). Here, we
would like to mention that our  result was already included in the
papers by Calabi when $m=1$ \cite{Ca},  by Mok \cite{Mo1} \cite{Mo2}
when $N_1=\cdots=N_m$, and by Ng \cite{Ng1} \cite{Ng3} when $m=2$
and $N_1, N_2<2n$.

As  in the work of Mok \cite{Mo1}, our proof of the theorem is also
based on the similar algebraic extension theorem derived in
\cite{Mo2} and Mok-Ng \cite{MN}. However, different from the case
considered in \cite{Mo1} \cite{Ng2}, the properness of a factor of
$F$ does not immediately imply the linearity of that factor; for the
classical linearity theorem does not hold anymore for proper
rational mappings from ${\BB^n}$ into ${\BB^N}$ with $N>2n-2$. (See
[Hu1]). Hence, the cancelation argument as in [Mok1] [Ng3] seems to
be difficult to apply in our setting.

 In our proof of Theorem 1.1, a major step is to
prove that a non-constant component $F_j$ of $F$ must be proper from $\BB^n$ into $\BB^{N_j}$,
 using the multi-valued holomorphic continuation technique.
This then reduces the proof of Theorem 1.1 to the case when all
components are proper. Unfortunately, due to the non-constancy for
the conformal factors $\lambda_j(z,\-{z})$ and $\lambda(z,\-{z})$,
it is not immediate that each component must also be conformal (and
thus $\lambda_j$ must be a constant multiple of $\lambda$) with
respect to the normalized Bergman metric. However, we observe that
the blowing-up rate for the Bergman metric of ${\BB^n}$ with $n\ge
2$ in the complex normal direction is twice of that along the
complex tangential direction, when approaching the boundary. From
this, we will be able to derive an equation regarding the CR
invariants associated to the map at the boundary of the ball.
Lastly, a linearity criterion of Huang in [Hu1] can be applied to
simultaneously conclude the linearity of all components.

 We mention that in the context of Corollary 1.2, namely, when each
conformal factor is assumed to be constant,
 the proof used to prove Theorem 1.1 can be  further simplified as told by Mok and Ng in their private
communications. In this case, one can work
 directly on the K\"ahler potential functions instead of the hyperbolic metrics.
 However, when the conformal factors are not constant, and $\partial\-{\partial}$-lemma cannot be applied and
  the metric equation (which can be regarded as differential equations on the map)
 does not lead to  the functional equation on the components of the map.
 We  appreciate very much many valuable comments of Mok and Ng to the earlier version of this paper,
 especially, for telling us how to essentially simplify the proof of a key lemma (Lemma 2.2) through the consideration of the
metric  potential functions. Their very helpful comments lead to the
present version.

\bigskip
{\bf Acknowledgement}: The  present work is written under the
guidance of  Professor Xiaojun Huang. The authors are deeply
indebted to Professor Huang, especially, for formulating the problem
and for suggesting the method  used in the paper. The authors also
would like to thank very much N. Mok and, in particular, S. Ng, for
many stimulating discussions and conversations which inspired the
present work as well as many valuable suggestions and comments. In
fact, this work was originated  by reading the work of  Mok ([Mo1])
and Ng ([Ng1-3]). The authors acknowledge the partial financial
support of NSF-0801056 for the summer research project (through
Professor Huang). Part of the work was done when the first author
was visiting Erwin Schr\"{o}dinger Institute in Vienna, Austria in the
fall semester of 2009. He also would like to thank members in the
institute, especially, Professor Lamel for the invitation and
hospitality. Finally, the authors would like to thank the referees for the valuable comments.
\bigskip

\section{ Bergman metric and  proper rational maps}
Let $\BB^n$ and $ds^2_n$  be the unit ball and its normalized
Bergman metric, respectively, as defined before.
Denote by ${\HH}^n\subset{\CC}^n$ the Siegel upper half space.
Namely, $\HH^n=\{(z,w)\in \CC^{n-1}\times \CC: \Im w- |z|^2>0\}$.
Here, for $m$-tuples $a, b$, we write dot product $a\cdot
b=\sum_{j=1}^{m}a_jb_j$ and  $|a|^2=a \cdot \bar a$. Recall  the
following  Cayley transformation

\begin {equation}\label{cayley}
\rho_n(z, w) = \bigg( \frac{2z}{1- i w},\ \frac{1+i w}{1- iw}
\bigg).
\end{equation}
 Then $\rho_n$ biholomorphically maps  ${\HH}^n$ to
${\BB}^n$, and biholomorphically maps $\p\HH^n$, the Heisenberg hypersurface, to $\p\BB^n
\backslash \{(0, 1) \}$. Applying the Cayley transformation, one can
 compute the normalized Bergman metric on $\HH^n$ as follows:
\begin{equation}\label{metric1}
\begin{split}
ds^2_{\HH^n}&=\sum_{j,k<n}\frac{\delta_{jk}(\Im w-|z|^2)+\bar z_j z_k}{(\Im
w-|z|^2)^2}dz_j\otimes d\bar z_k+\frac{dw\otimes d\bar
w}{4(\Im w-|z|^2)^2}\\
&\ \ \ \ +\sum_{j<n}\frac{\bar z_jdz_j\otimes d\bar w}{2i(\Im
w-|z|^2)^2}-\sum_{j<n}\frac{z_j dw\otimes d\bar z_j}{2i(\Im
w-|z|^2)^2}.
\end{split}
\end{equation}
One can easily check that
\begin{align*}
&L_j=\frac{\p}{\p z_j}+2i\bar z_j\frac{\p}{\p w}, j=1,\ldots,n-1.\\
&\overline{L_j}=\frac{\p}{\p \bar z_j}-2i z_j\frac{\p}{\p \bar{w}}, j=1,\ldots,n-1.\\
&T=2(\frac{\p}{\p w}+\frac{\p}{\p \bar w})
\end{align*}
span the complexified tangent vector bundle of $\p\HH^n$. (See, for
instance, \cite{BER}, \cite{Hu2}, \cite{Hu3} [HJX] [JX] [Mi].)

Let $F$ be a rational proper holomorphic map from $\HH^n$ to
$\HH^N$. By a result of Cima-Suffridge [CS], $F$ is holomorphic in a
neighborhood of $\p{\HH^n}$. Assign the weight of $w$ to be 2 and
that of $z$ to be 1. Denote by $o_{wt}(k)$ terms with weighted degree
higher than $k$ and by $P^{(k)}$ a function of weighted degree $k$.
For $p_0=(z_0,w_0)\in \p{\HH^n}$, write
$\sigma^{0}_{p_0}:(z,w)\rightarrow (z+z_0,w+w_0+2i z\cdot \bar z_0)$
for the standard Heisenberg translation. The following normalization
lemma  will be used here:

\begin{lemma}\label{moser}
[Hu2-3] For any $p\in \p{\HH^n}$, there is an element $\tau\in
Aut(\HH^{N+1})$ such that the map
$F^{**}_{p}=((f^{**}_p)_1(z),\cdots,(f^{**}_p)_{n-1}(z),\phi_p^{**},
g_p^{**})= (f_p^{**},\phi_p^{**}, g_p^{**})=\tau\circ F\circ
\sigma_p^0$ takes the following normal form:
\begin{align*}
&f^{**}_p(z,w)=z+ \frac{i}{2}a^{(1)}(z)w+o_{wt}(3),\\
&\phi^{**}_p(z,w)=\phi^{(2)}(z)+o_{wt}(2),\\
&g^{**}_p(z,w)=w+ o_{wt}(4)
\end{align*}
with
\begin{equation}\label{moser2}(\bar z\cdot
a^{(1)}(z))|z|^2=|\phi^{(2)}(z)|^2.
\end{equation}
\end{lemma}

In particular, write
$(f^{**}_p)_l(z)=z_j+\frac{i}{2}\sum_{k=1}^{n-1}a_{lk}z_kw+o_{wt}(3).$
Then, $(a_{lk})_{1\le l,k\le n-1}$ is an $(n-1)\times(n-1)$
semi-positive Hermitian metrix. We next present the following key
lemma for our proof of Theorem 1.1:



\begin{lemma}\label{smooth-ball}
Let $F$ be a proper rational map from $\BB^n$ to $\BB^N$. Then
\begin{equation}\label{X}
X:=ds^2_{n}-F^*(ds^2_{N}),
\end{equation}
is a semi-positive real analytic symmetric (1,1)-tensor  over
$\BB^n$ that  extends also to a real analytic  (1,1)-tensor in a
small neighborhood of  $\p\BB^n$ in $\CC^n$.
\end{lemma}

\medskip

{\noindent \bf Proof of lemma \ref{smooth-ball}}: Our original proof was
largely simplified by Ng \cite{Ng4} and Mok [Mo4] by considering the
potential $-\log(1-\left\| F(z) \right\|^2)$ of the pull-back metric
$F^*({ds^2_N})$ as follows:
 Since $1-\left\| F(z) \right\|^2$ vanishes identically
on $\p\BB^n$ and since $1-\left\| z \right\|^2$ is a defining equation for $\p\BB^n$, one
obtains

$$1-\left\| F(z) \right\|^2 = (1-\left\| z \right\|^2) \varphi(z)$$ for a real analytic function $\varphi(z)$.

Since $\rho:=\left\| F(z) \right\|^2 -1$ is subharmonic over
${\BB^n}$ and has maximum value $0$ on the boundary, applying the
classical Hopf lemma, we conclude that $\varphi(z)$ cannot vanish
at any boundary point of $\BB^n$. Apparently, $\varphi(z)$ cannot
vanish inside $\BB^n$. Therefore, $X=\sqrt{-1}\p\bar\p \log
\varphi(z)$ is real analytic on an open neighborhood of
$\overline{\BB^n}$. The semi-positivity of $X$ over ${\BB^n}$ is an
easy consequence of the Schwarz lemma.
\endpf

\medskip

Applying the Cayley transformation (and also a rotation
transformation when handling the regularity near  $(0,1)$), we have
the following corollary:

\begin{corollary}\label{smooth}
Let $F$ be a rational proper holomorphic map from $\HH^n$ to
$\HH^N$. Then
\begin{equation}\label{X}
X:=ds^2_{\HH^n}-F^*(ds^2_{\HH^N}),
\end{equation}
is a semi-positive real analytic symmetric (1,1)-tensor  over
$\HH^n$ that  extends also to a real analytic  (1,1)-tensor in a
small neighborhood of $\p\HH^n$ in $\CC^n$.
\end{corollary}

\medskip
The boundary value of $X$ is an intrinsic CR invariant associated
with the equivalence class of the map $F$. Next, we compute $X$ in
the normal coordinates at the boundary point.

Write $t=\Im w-|z|^2$ and $H=\Im g- |\tilde f|^2$. Here $(\tilde f,
g)$ denotes the map between Heisenberg hypersurfaces. Write $o(k)$
for terms whose degrees with respect to $t$ are higher than $k$. For
a real analytic function $h$ in $(z,w)$, we use $h_{z}, h_{w}$ to
denote the derivatives of $h$ with respect to $z, w$. By replacing
$w$ by $u+i(t+|z|^2)$, $H$ can also be regarded as an analytic
function on $z,\bar z, u, t$. The following lemma gives an
asymptotic behavior of $H$ with respect to $t$:

\medskip
\begin{lemma}\label{P1P2}
 $H(z,\-{z},u,t)= (g_w-2i\w f_w \cdot \- {\w f} )|_{t=0}t -(2|\w
f_{w}|^2)|_{t=0}t^2+\frac{1}{3}(-\frac{1}{2}g_{w^3}+3i \w f_w\cdot
\overline{\w f_{w^2}}+i\w f_{w^3}\cdot \- {\w f})|_{t=0}t^3+o(3)$.
\end{lemma}
\medskip
{\noindent \bf Proof of Lemma \ref{P1P2}}: Write $H=H(z, \-z,
u+i(t+|z|^2), u-i(t+|z|^2))$. Since $F$ is proper, $H$, as a
function of $t$ with parameters  $\{z,u\}$, can be written as $P_1t
+ P_2 t^2+ P_3 t^3+o(3)$, where $P_1, P_2, P_3$ are analytic in
$(z,\-z, u)$. Then

\begin{equation}\label{P1}
\begin{split}
 P_1&=\frac{\p H(z, \-z,
u+i(t+|z|^2), u-i(t+|z|^2))}{\p t}\bigg|_{t=0}\\
&=iH_w-i  H_{\-w}\bigg|_{t=0}\\
&=\frac{1}{2}(g_w+\overline{ g_w})+i(\w f\cdot \overline {\w
f_w}-\bar {\w f}\cdot \w f_w)\bigg|_{t=0},
\end{split}
\end{equation}

\begin{equation}\label{P2}
\begin{split}
 P_2&= \frac{1}{2}\frac{\p^2 H(z, \-z,
u+i(t+|z|^2), u-i(t+|z|^2))}{\p t^2}\bigg|_{t=0}\\
&=\frac{1}{2}(-H_{w^2}+2H_{w\bar w}-H_{\- w^2})\bigg|_{t=0}\\
&=\frac{1}{2}(\frac{i}{2}g_{w^2}-\frac{i}{2}\overline{g_{w^2}}-2|\tilde
f_w|^2+\tilde f_{w^2}\cdot \bar{\w f}+\w f\cdot \overline{\w
f_{w^2}})\bigg|_{t=0},
\end{split}
\end{equation}
and
\begin{equation}\label{P3}
\begin{split}
P_3&=\frac{1}{6}\frac{\p^3H(z,\bar z, u+i(t+|z|^2), u-i(t+|z|^2))}{\p t^3}\bigg|_{t=0}\\
&=\frac{1}{6}(-iH_{w^3}+3iH_{w^2\-w}-3iH_{\bar w^2w}+iH_{\- w^3})\bigg|_{t=0}\\
&=\frac{1}{6}(-\frac{1}{2}g_{w^3}-\frac{1}{2}\overline{g_{w^3}}+i\w
f_{w^3}\cdot \- {\w f}-i \w f\cdot \overline{\w f_{w^3}}-3i \w
f_{w^2}\cdot \overline{\w f_{w}}+3i \w f_w\cdot \overline{\w
f_{w^2}})\bigg|_{t=0}.
\end{split}
\end{equation}
On the other hand, applying $T, T^2, T^3$
to the defining equation $g-\bar g= 2i \w f\cdot \bar {\w f}$, we have
\begin{align}
&g_w-\overline{ g_w}-2i(\w f_w\cdot \bar {\w f}+\w f\cdot
\overline{\w f_w})=0,\label{T}\\
&g_{w^2}-\overline{g_{w^2}}-2i(\w f_{w^2}\cdot \-{\w
f}+\overline{\w f_{w^2}}\cdot \w f+2|\w f_w|^2)=0, \label{TT}\\
&g_{w^3}-\overline{g_{w^3}}-2i(\w f_{w^3}\cdot \w f+ \overline{\w
f_{w^3}}\cdot \w f+3\w f_{w^2}\cdot \overline {\w f_w}+3 \w
f_{w}\cdot \overline{\w f_{w^2}})=0.\label{TTT}
\end{align}
over $\Im w=|z|^2$.

 Substituting  (\ref{T}), (\ref{TT}) and (\ref{TTT}) into (\ref{P1}), (\ref{P2}) and (\ref{P3}), we get
\begin{equation}\label{formula}
\begin{split}
&P_1=g_w-2i\w f_w\cdot \bar {\w f}\bigg|_{t=0},\\
&P_2= -2|\w f_{w}|^2\bigg|_{t=0},\\
&P_3=\frac{1}{3}(-\frac{1}{2}g_{w^3}+3i \w f_w\cdot \overline{\w
f_{w^2}}+i\w f_{w^3}\cdot \- {\w f})\bigg|_{t=0}. \hspace{3cm}
\endpf
\end{split}
\end{equation}

We remark that by the Hopf Lemma, it follows easily that $P_1\neq 0$
along $\p{\HH^n}$.

\bigskip
We next write $X=X_{jk}dz_j\otimes d\bar z_k+X_{jn}dz_j\otimes d\bar
w+ X_{nj}dw\otimes d\bar z_j+X_{nn}dw\otimes d\bar w$. By making use
of Lemma \ref{moser}, we shall compute in the next proposition the
values of $X$ at the origin. The proposition might be of independent
interest, as the CR invariants in the study of proper holomorphic
maps between Siegel upper half spaces are related to the CR geometry
of the map.

\begin{proposition}\label{zero}
Assume that $F=(\tilde f, g)=(f_1,\ldots, f_{N-1}, g):\HH^n\rightarrow \HH^N$
is a proper rational holomorphic map that satisfies the normalization (at the
origin) stated in Lemma \ref{moser}. 
Then
\begin{align*}
&X_{jk}(0)=-2i(f_k)_{z_jw}(0)=a_{kj},\\
&X_{jn}(0)=\overline{X_{nj}}(0)=\frac{3i}{4}\overline{(f_j)_{w^2}}(0)+\frac{1}{8}g_{z_jw^2}(0), \\
&X_{nn}(0)=\frac{1}{6}g_{w^3}(0).
\end{align*}
\end{proposition}
\medskip
\noindent{\bf Proof of Proposition \ref{zero}:}\  Along the
direction of $dz_j\otimes d \bar z_k$, collecting the coefficient of
$t^2$ in the Taylor expansion of $H^2X$ with respect to $t$, we get
\begin{align*}
P_1^2X_{jk}(0)=&\bigg[(2P_1P_2\delta_{jk}+(P_2^2+2P_1P_3)\bar z_j
z_k)-\frac{1}{2}\big{\{}2iP_1(\w f_{wz_{j}}\cdot \overline{\w f_{z_{k}}}-\w f_{z_{j}}\cdot
\overline{\w f_{w{z_{k}}}}) +2P_2 \w f_{z_{j}}\cdot \overline{\w f_{z_{k}}}\\
&-(\overline{\w f_{w^2}}\cdot \w f_{z_{j}})(\w f\cdot \overline{\w f_{z_{k}}})+2
(\overline{\w f_w}\cdot \w f_{z_{j}})(\w f_w\cdot \overline{\w f_{z_{k}}})+2
(\overline{\w f_w}\cdot \w f_{w{z_{j}}})(\w f\cdot \overline{\w f_{z_{k}}})\\
&-2
(\overline{\w f_w}\cdot \w f_{z_{j}})(\w f\cdot \overline{\w
f_{{z_{k}}w}})-(\-{\w f}\cdot\w f_{z_{j}})(\w f_{w^2}\cdot \overline{\w f_{z_{k}}})-2(\-{\w
f}\cdot \w f_{{z_{j}}w})(\w f_w\cdot \overline{\w f_{z_{k}}})\\
& +2 (\-{\w f}\cdot
\w f_{z_{j}})(\w f_w\cdot \overline{\w f_{{z_{k}}w}})-(\w
f_{{z_{j}}w^2}\cdot\-{ \w f})(\w f\cdot \overline{\w f_{z_{k}}})+2(\-{\w f}\cdot \w f_{{z_{j}}w})(\w f\cdot \overline{\w f_{{z_{k}}w}}) \\
&-(\-{\w
f}\cdot \w f_{z_{j}})(\w f\cdot \overline{\w f_{{z_{k}}w^2}})
-\frac{1}{4}g_{{z_{j}}w^2}\overline{g_{z_{k}}}+\frac{1}{2}g_{{z_{j}}w}\overline{g_{{z_{k}}w}}-\frac{1}{4}g_{z_{j}}\overline{g_{{z_{k}}w^2}}
+\frac{i}{2}(\overline{\w f_{w^2}}\cdot \w f_{z_{j}})\bar
g_{z_{k}}\\&-i(\overline{\w f_w}\cdot\w f_{w{z_{j}}})\overline{g_{z_{k}}}+i(\overline{\w
f_w}\cdot \w f_{z_{j}})\overline{g_{{z_{k}}w}}+\frac{i}{2}(\-{\w f}\cdot \w
f_{{z_{j}}w^2})\overline{g_{z_{k}}} -i(\-{\w f}\cdot\w
f_{{z_{j}}w})\overline{g_{{z_{k}}w}}\\
&+\frac{i}{2}(\-{\w f}\cdot \w
f_{z_{j}})\overline{g_{{z_{k}}w^2}}-\frac{i}{2}(\w f_{w^2}\cdot
\overline{\w f_{z_{k}}})g_{z_{j}}-i(\w f_w\cdot \overline{\w f_{{z_{k}}}})g_{{z_{j}}w}+i(\w f_w\cdot
\overline{\w f_{{z_{k}}w}})g_{z_{j}}\\
&-\frac{i}{2}(\w f\cdot \overline{\w
f_{z_{k}}})g_{{z_{j}}w^2}+i (\w f\cdot\overline{\w f_{{z_{k}}w}})g_{{z_{j}}w}-\frac{i}{2}(\w
f\cdot \overline{\w f_{{z_{k}}w^2}})g_{z_{j}} \big{\}}\bigg]\bigg|_{t=0}.
\end{align*}
Letting $(z,w)=0$ and applying the normalization condition as stated
in Lemma \ref{moser}, we  have
\begin{align*}
X_{jk}(0)=\frac{\p a^{(1)}_k(z)}{\p z_j}=a_{kj}.
\end{align*}
Similarly, considering the coefficients of $t^2$ along $dz_j\otimes
d\-w$ and $dw\otimes d\-w$, respectively,  we have

\begin{align*}
 P_1^2X_{jn}(0)=&\bigg[(-iP_1P_3-\frac{i}{2}P_2^2)\bar z_j
-\frac{1}{2}\{2iP_1(\w f_{w{z_{j}}}\cdot \overline{\w f_w}-\w f_{z_{j}}\cdot
\overline{\w f_{w^2}}) +2P_2 \w f_{z_{j}}\cdot \overline{\w f_w}\\
&-(\overline{\w f_{w^2}}\cdot \w f_{z_{j}})(\w f\cdot \overline{\w f_w})+2
(\overline{\w f_w}\cdot \w f_{z_{j}})(\w f_w\cdot \overline{\w f_w})+2
(\overline{\w f_w}\cdot \w f_{w{z_{j}}})(\w f\cdot \overline{\w f_w})\\
&-2
(\overline{\w f_w}\cdot \w f_{z_{j}})(\w f\cdot \overline{\w
f_{w^2}})-(\-{\w f}\cdot\w f_{z_{j}})(\w f_{w^2}\cdot \overline{\w f_w})-2(\-{\w
f}\cdot \w f_{{z_{j}}w})(\w f_w\cdot \overline{\w f_w})
\end{align*}
\begin{align*}
& +2 (\-{\w f}\cdot
\w f_{z_{j}})(\w f_w\cdot \overline{\w f_{w^2}})-(\w
f_{{z_{j}}w^2}\cdot\-{ \w f})(\w f\cdot \overline{\w f_w})+2(\-{\w f}\cdot \w f_{{z_{j}}w})(\w f\cdot \overline{\w f_{w^2}})
-(\-{\w f}\cdot \w f_{z_{j}})(\w f\cdot \overline{\w f_{w^3}})\\
&-\frac{1}{4}g_{{z_{j}}w^2}\overline{g_w}+\frac{1}{2}g_{{z_{j}}w}\overline{g_{w^2}}-\frac{1}{4}g_{z_{j}}\overline{g_{w^3}}
+\frac{i}{2}(\overline{\w f_{w^2}}\cdot \w f_{z_{j}})\overline{g_w}-i(\overline{\w f_w}\cdot\w f_{w{z_{j}}})\overline{g_w}+i(\overline{\w
f_w}\cdot \w f_{z_{j}})\overline{g_{w^2}}\\
&+\frac{i}{2}(\-{\w f}\cdot \w
f_{{z_{j}}w^2})\overline{g_w} -i(\-{\w f}\cdot\w
f_{{z_{j}}w})\overline{g_{w^2}}+\frac{i}{2}(\-{\w f}\cdot \w
f_{z_{j}})\overline{g_{w^3}}-\frac{i}{2}(\w f_{w^2}\cdot
\overline{\w f_w})g_{z_{j}}-i(\w f_w\cdot \overline{\w f_{w}})g_{{z_{j}}w}\\
&+i(\w f_w\cdot
\overline{\w f_{w^2}})g_{z_{j}}-\frac{i}{2}(\w f\cdot \overline{\w
f_w})g_{{z_{j}}w^2}+i (\w f\cdot\overline{\w
f_{w^2}})g_{{z_{j}}w}-\frac{i}{2}(\w f\cdot \overline{\w f_{w^3}})g_{z_{j}}
\}\bigg]\bigg|_{t=0}
\end{align*}
and

\begin{align*}
 P_1^2X_{nn}(0)=&\bigg[\frac{1}{4}(2P_1P_3+P_2^2)-\frac{1}{2}\{2iP_1(\w f_{w^2}\cdot \overline{\w f_w}-\w f_w\cdot
\overline{\w f_{w^2}}) +2P_2 \w f_w\cdot \overline{\w f_w}-(\overline{\w f_{w^2}}\cdot \w f_w)(\w f\cdot \overline{\w f_w})\\
&+2
(\overline{\w f_w}\cdot \w f_w)(\w f_w\cdot \overline{\w f_w})+2
(\overline{\w f_w}\cdot \w f_{w^2})(\w f\cdot \overline{\w f_w})-2
(\overline{\w f_w}\cdot \w f_w)(\w f\cdot \overline{\w
f_{w^2}})\\
&-(\-{\w f}\cdot\w f_w)(\w f_{w^2}\cdot \overline{\w f_w})-2(\-{\w
f}\cdot \w f_{w^2})(\w f_w\cdot \overline{\w f_w}) +2 (\-{\w f}\cdot
\w f_w)(\w f_w\cdot \overline{\w f_{w^2}})-(\w
f_{w^3}\cdot\-{ \w f})(\w f\cdot \overline{\w f_w})\\
&+2(\-{\w f}\cdot \w f_{w^2})(\w f\cdot \overline{\w f_{w^2}})
-(\-{\w f}\cdot \w f_w)(\w f\cdot \overline{\w f_{w^3}})
-\frac{1}{4}g_{w^3}\overline{g_w}+\frac{1}{2}g_{w^2}\overline{g_{w^2}}-\frac{1}{4}g_w\overline{g_{w^3}}
\\
&+\frac{i}{2}(\overline{\w f_{w^2}}\cdot \w f_w)\bar
g_w-i(\overline{\w f_w}\cdot\w f_{w^2})\overline{g_w}+i(\overline{\w
f_w}\cdot \w f_w)\overline{g_{w^2}}+\frac{i}{2}(\-{\w f}\cdot \w
f_{w^3})\overline{g_w} -i(\-{\w f}\cdot\w
f_{w^2})\overline{g_{w^2}}\\
&+\frac{i}{2}(\-{\w f}\cdot \w
f_w)\overline{g_{w^3}}-\frac{i}{2}(\w f_{w^2}\cdot
\overline{\w f_w})g_w-i(\w f_w\cdot \overline{\w f_{w}})g_{w^2}\\
&+i(\w f_w\cdot
\overline{\w f_{w^2}})g_w-\frac{i}{2}(\w f\cdot \overline{\w
f_w})g_{w^3}+i (\w f\cdot\overline{\w
f_{w^2}})g_{w^2}-\frac{i}{2}(\w f\cdot \overline{\w f_{w^3}})g_w
\}\bigg]\bigg|_{t=0}.
\end{align*}

Let $(z,w)=0$. It follows that
\begin{align*}
&X_{jn}(0)=\frac{3i}{4}\overline{(f_j)_{w^2}}(0)+\frac{1}{8}g_{{z_{j}}w^2}(0), \\
&X_{nn}(0)=\frac{1}{6}g_{w^3}(0),
\end{align*}
for $g_{w^3}(0)=\overline{g_{w^3}(0)}$ by (\ref{TTT}).\ \ \ \ \ \ \endpf

\bigskip

Making use of the computation in Proposition 2.5, we give a proof of
Theorem \ref{main1} in the case when each component extends as a
proper holomorphic map. Indeed, we prove a slightly more general
result than what we will need later as follows.
\bigskip

\begin{proposition}\label{mainprop}
Let
\begin{align*} &F=(F_1,\ldots,
F_m): \BB^n
\rightarrow \BB^{N_1}\times \cdots \times
\BB^{N_m}
\end{align*}
be a holomorphic isometric embedding up to conformal factors
$\{\lambda(z,\-{z});\lambda_1(z,\-{z}),\cdots, \lambda_m(z,\-{z})\}$
in the sense that
 $$\lambda(z,\-{z}) ds^2_n=\sum_{j=1}^{m}\lambda_j(z,\overline{z})
F_j^*(ds^2_{N_j}).$$ Here for each $j$, $\lambda(z, z), \lambda_j(z,
\overline{z})$ are positive 
$C^2$-smooth functions
over $\overline{\BB^n}$, and $F_j$ is a proper rational map from
$\BB^n$ into $\BB^{N_j}$ for each $j$. Then $\lambda(z,
\-{z})\equiv\sum_{j=1}^{m} \lambda_j(z, \overline{z})$ over
$\overline{\BB^n}$, and for any $j$, $F_j$ is a totally geodesic
embedding from $\BB^n$ to $\BB^{N_j}$.
\end{proposition}
\medskip
\noindent{\bf Proof of Proposition \ref{mainprop}:} \ After applying
the Cayley transformation and considering
$\left((\rho_{N_1})^{-1},\cdots,(\rho_{N_m})^{-1}\right)\circ F\circ
\rho_n$ instead of $F$, we can assume, without loss of generality,
that  $$F=(F_1,\ldots, F_m): \HH^n
\rightarrow
\HH^{N_1}\times \cdots \times \HH^{N_m}
$$ is an isometric map up to conformal factors
$\{\lambda(Z,\-{Z});\lambda_1(Z,\-{Z}),\cdots,
\lambda_m(Z,\-{Z})\}$ in the sense that
 $$\lambda(Z,\-{Z}) ds^2_{\HH^n}=\sum_{j=1}^{m}\lambda_j(Z,\overline{Z})
F_j^*(ds^2_{\HH^{N_j}}).$$
Also, each $F_j$ is a proper rational map from $\HH^n$ into
$\HH^{N_j}$, respectively. Here we write $Z=(z,w)$. Moreover, we can
assume, without loss of generality, that each component $F_j$ of $F$
satisfies the normalization condition as in Lemma 2.1. Since $F$ is an isometry, we have
\begin{equation}\label{0000}
\lambda(Z, \bar Z) ds^2_{\HH^n}=\sum_{j=1}^{m}\lambda_j(Z, \bar Z)
F^*_j(ds^2_{\HH^{N_j}}),\ \hbox{or } (\lambda(Z, \bar
Z)-\sum_{j=1}^{m}\lambda_j(Z, \bar Z))ds^2_{\HH^n}+
\sum_{j=1}^{m}\lambda_j(Z, \bar Z) X(F_j)=0.
\end{equation}
Here, we write $X(F_j)=ds^2_{\HH^n}-F^*_j(ds^2_{\HH^{N_j}})$.
Collecting the coefficient of $dw\otimes d\bar w$, one has
\begin{equation}\label{00000}
\frac{\lambda(Z, \bar Z)-\sum_{j=1}^{m} \lambda_j(Z, \bar Z)}{4(\Im
w-|z|^2)^2} + \sum_{j=1}^{m}\lambda_j(Z,\-{Z})
(X(F_j))_{nn}=0.\end{equation}
 Since $X(F_j)$ is smooth up to  $\p
{\HH^n}$, we see that $\lambda(Z, \bar Z)-\sum_{j=1}^{m}
\lambda_j(Z, \bar Z)=O(t^2)$ as $Z=(z,w)(\in {\HH^n})\ra 0$, where
$t=\Im w-|z|^2$.
However, since the $dz_l\otimes d\-{z_k}$-component of
$ds^2_{\HH^n}$ blows up at the rate of $o(\frac{1}{t^2})$ as
$(z,w)(\in{\HH^n})\ra 0$, collecting the coefficients of the
$dz_l\otimes d \bar z_k$-component in (\ref{0000}) and then letting
$(z,w)(\in {\HH^n})\ra 0$, we conclude that, for any $1\le l, k\le
n-1$,
$$\sum_{j=1}^m \lambda_j(0, 0) (X(F_j))_{k{l}}(0)=0.$$ By
Proposition 2.5, we have $\sum_{j=1}^m \lambda_j(0, 0)
a^j_{l{k}}(0)=0,$ where $a^j_{kl}$ is associated with $F_j$ in the
expansion of $F_j$ at $0$ as in Lemma 2.1.  Since $(a^j_{kl})_{1\le l, k\le n-1}$ is a semi-positive matrix and
$\lambda_j(0, 0)>0$, it follows immediately that
$a_{lk}^j(0)=0$ for all $j, k,l$. Namely, $F_j=(z,w)+O_{wt}(3)$ for
each $j$.


Next, for each $p\in \p \HH^n$, let $\tau_j \in Aut(\HH^{N})$ be
such that $(F_j)_{p}^{**}=\tau_j \circ F_j \circ \sigma_p^0$ has the
normalization  as in Lemma \ref{moser}. Let
$\tau=(\tau_1,\cdots,\tau_m)$. Note that
$F_p^{**}:=((F_1)_p^{**},\cdots,(F_m)_p^{**})=\tau\circ F\circ
\sigma_p^0$ is still an isometric map satisfing the condition as in
the proposition. Applying the just presented argument to $F_p^{**}$,
we conclude that $(F_j)_p^{**}=(z, 0, w)+O_{wt}(3).$ By Theorem 4.2
of [Hu2], this implies that $F_j=(z, 0, w)$. Namely, $F_j$ is a totally
geodesic embedding. In particular, we have $X(F_j) \equiv 0$. This
also implies that $\lambda\equiv \sum_{j=1}^{m} \lambda_j$ over
$\BB^n$. The proof of Proposition \ref{mainprop} is complete.
$\endpf$

\section{Proof of Theorem \ref{main1}}
In this section, we give a proof of Theorem \ref{main1}. As in the
theorem, we let $U\subset \BB^n$ be a connected  open subset. Let
\begin{align*}
&F=(F_1,\ldots, F_m): U
\rightarrow \BB^{N_1}\times
\cdots \times \BB^{N_m}
\end{align*}
be a holomorphic isometric embedding up to conformal factors
$\{\lambda(z,\-{z});\lambda_1(z,\-{z}),\cdots, \lambda_m(z,\-{z})\}$
in the sense that
 $$\lambda(z,\-{z}) ds^2_n=\sum_{j=1}^{m}\lambda_j(z,\overline{z})
F_j^*(ds^2_{N_j}).$$
 Here $\lambda_j(z, \bar z),
\lambda(z, \bar z)>0$ are smooth Nash algebraic functions;
$ds_{n}^2$ and $ds_{N_j}^2$ are the Bergman metrics of $\BB^n$ and
$\BB^{N_j}$, respectively; and $F_j$ is a holomorphic map from $U$
into $\BB^{N_j}$ for each $j$. For the proof of Theorem \ref{main1},
we can assume without loss of generality that none of the $F_j$'s is
a constant map. Following the idea in \cite{MN}, we can show that
$F$ extends to an algebraic map over $\CC^n$. (For the convenience
of the reader, we include the detailed argument in the appendix.)
Namely, for each (non-constant) component $f_{j,l}$ of $F_j$, there
is an irreducible polynomial $P_{j,l}(z,
X)=a_{j,l}(z)X^{m_{jl}}+\ldots$ in $(z, X)\in \CC^n\times \CC$ of
degree $m_{jl}\ge 1$ in $X$ such that $P_{j,l}(z, f_{j,l})\equiv 0$
for $z\in U$.

We will proceed to show that, for each $j$, {$F_j$ extends to a
proper rational map from $\BB^n$ into $\BB^{N_j}$.}  For this
purpose, we let $R_{j,l}(z)$ be the resultant of $P_{j,l}$ in
$X$ and let
 $E_{j,l}=\big\{R_{j,l}\equiv 0, a_{j,l}\equiv 0\big\}$, $E=\cup E_{j,l}$.
 Then $E$ defines a proper affine-algebraic subvariety of $\CC^n$.
 For any continuous curve $\gamma: [0, 1]\rightarrow \CC^n\setminus E$ where $\gamma(0)\in U$, $F$
 can be continued holomorphically along $\gamma$ to get a germ of holomorphic map at $\gamma(1)$. Also,
 if $\gamma_1$ is homotopic to $\gamma_2$ in $\CC^n\setminus E$, $\gamma_1(0)=\gamma_2(0)\in U$ and $\gamma_1(1)=\gamma_2(1)$,
 then continuations of $F$ along $\gamma_1$ and $\gamma_2$ are the same at $\gamma_1(1)=\gamma_2(1)$.
  Now let $p_0\in U$ and $p_1\in \partial \BB^n\setminus E$. Let $\gamma(t)$ be a smooth simple curve connecting $p_0$ to $p_1$
  and $\gamma(t)\notin \partial \BB^n$ for $t\in (0,1)$. Then each $F_j$ defines a holomorphic map in a connected neighborhood
   $V_{\gamma}$ of $\gamma$ by continuing along $\gamma$ the initial germ of $F_j$ at $p_0$.
   (We can also assume that $V_\gamma\cap {\BB^n}$ is connected.)  Let
 \begin{align*}
 S_{\gamma}=\big\{p\in V_{\gamma}:\ \|F_j(p)\|=1 \text{\ for some} \ j\}.
 \end{align*}
Then $S_{\gamma}$ is a real analytic (proper) subvariety of
$V_{\gamma}$. We first prove

\begin{claim}\label{dimension} When $V_\gamma$ is sufficiently close
to $\gamma$, $dim(S_{\gamma}\cap \BB^n)\le 2n-2$.
\end{claim}
{\noindent \bf Proof of Claim \ref{dimension}:}\ Supposing otherwise we are going to deduce a contradiction.
Assume that $t_0\in (0,1]$ is the first point such that for a
certain $j$, the local variety defined by $\|F_{j}(z)\|^2=1$ near
$p^*=\gamma(t_0)$ has real dimension $2n-1$ at $p^*$. Let $\Sigma_0$ be an irreducible component of the germ of the real analytic subvariety $S_\gamma$ at $p^*$ of real codimension 1 and let $\Sigma$ be a connected locally closed subvariety of $\mathbb{B}^n$ representing the germ $\Sigma_0$ at $p^*$. Since any real
analytic subset of real codimension two inside a connected open set does
not affect the connectivity, by slightly changing $\gamma$ without
changing its homotopy type and terminal point, we can assume that
$\gamma(t)\not\in S_\gamma$ for any $t<t_0$. Hence, $p^*$  also lies
on the boundary of the connected component $\hat V$ of
$(V_{\gamma}\cap \BB^n)\setminus S_{\gamma}$, that contains
$\gamma(t)$ for $t<t_0$ and $\Sigma$ also 
lies in the boundary of
$\hat V$. Now, for any $p\in \Sigma$, let $q(\in \hat V)\rightarrow
p$, we have along $\{q\}$,
\begin{align*}
\lambda(z, \bar z) ds^2_n=\sum_j \lambda_j(z, \bar z) F_j^*(ds^2_{N_j}).
\end{align*}
Suppose that $j^{\sharp}$ is such that $\|F_{j^{\sharp}}(p)\|=1$ and
$\|F_j(z)\|<1$ for any $j$, $p\in \Sigma$ and $z\in \hat{V}$. Since
$p\in\Sigma\subset \BB^n$, $ds^2_n$ is a smooth Hermitian metric in an open neighborhood of $p$. For any $v=(v_1,\ldots, v_{\xi},\ldots, v_n) \in \mathbb{C}^n$ with $\|v\|=1$, 
it follows that
\begin{align*}
\left|\overline \lim_{q\rightarrow p}
F^*_{j^{\sharp}}(ds^2_{N_{j^{\sharp}}})(v, v)(q)\right| <\infty.
\end{align*}
On the other hand,
\begin{align*}
F_{j^{\sharp}}^*(ds^2_{N_{j^{\sharp}}})=\frac{\sum_{l,k}
\{\delta_{lk}(1-\|F_{j^{\sharp}}\|^2)+\bar f_{j^{\sharp},l}f_{j^{\sharp},k}\}d
f_{j^{\sharp},l}\otimes d\overline f_{j^{\sharp},k}}{(1-\|F_{j^{\sharp}}\|^2)^2}.
\end{align*}
It follows that 
\begin{align}\label{F*}
F_{j^{\sharp}}^*(ds^2_{N_{j^{\sharp}}})(v, v)(q)=\frac{\|\sum_{\xi}\frac{\partial
f_{j^{\sharp},l}}{\partial
z_{\xi}}(q)v_{\xi}\|^2}{1-\|F_{j^{\sharp}}(q)\|^2}+\frac{|\sum_{l,\xi}\overline
{f_{j^{\sharp},l}(q)}\frac{\p f_{j^{\sharp},l}}{\p
z_{\xi}}(q)v_{\xi}|^2}{(1-\|F_{j^{\sharp}}(q)\|^2)^2}.
\end{align}
Letting $q\rightarrow p$, since $1-\|F_{j^{\sharp}}(q)\|^2\rightarrow
0^+$, we get
\begin{align*}
\big \|\sum_{\xi}\frac{\p f_{j^{\sharp},l}(p)}{\p z_{\xi}}v_{\xi}\big\|^2=0.
\end{align*}
Thus
\begin{align*}
\frac{\p f_{j^{\sharp},l}(p)}{\p z_{\xi}}=0, \ \ \text{for} \ \ l=1,\ldots,
N_{j^{\sharp}}.
\end{align*}
Hence, we see  $dF_{j^{\sharp}}=0$ in a certain open subset of $\Sigma$.
Since $\Sigma$ is of real codimension 1 in $\mathbb{B}^n$, any non-empty open subset of $\Sigma$ is a uniqueness set for
holomorphic functions. Hence $F_{j^{\sharp}}\equiv const$. This is a
contradiction. $\endpf$

\medskip
Now, since $dim (S_{\gamma}\cap \BB^n) \le 2n-2$, we can always
slightly change $\gamma$ without changing the homotopy  type of
$\gamma$ in $V_{\gamma}\setminus E$ and the end point of $\gamma$ so
that $\gamma(t)\notin S_{\gamma}$ for any $t\in (0,1)$. Since
$\lambda(z, \bar z) ds^2_n=\sum^m_{j=1} \lambda_j(z, \bar z)
F_j^*(ds^2_{N_j})$ in $(V_{\gamma}\cap {\BB^n}) \setminus
S_{\gamma}$ and since $ds^2_n$ blows up when $q\in V_{\gamma}\cap
\BB^n$ approaches to $\p \BB^n$, we see that for each $q\in
V_{\gamma}\cap
\partial \BB^n$, $\|F_{j_q}(q)\|=1$ for some $j_q$. Hence, we can
assume without loss of generality, that there is a $j_0\ge 1$ such
that each of $F_1, \ldots, F_{j_0}$ maps a certain open piece
of $\p {\BB^n}$ into $\p \BB^{N_1},\ldots, \p \BB^{N_{j_0}}$, but
for $j>j_0$,
\begin{align*}
dim\{q\in \p \BB^n\cap V_{\gamma}: \|F_j(q)\|=1\}\le 2n-2.
\end{align*}
It follows from the Hopf lemma that $N_j\ge n$ for $j\le j_0$. By the
results of Forstneric [Fo] and Cima-Suffridge [CS], $F_j$ extends to
a rational proper holomorphic map from $\BB^n$ into $\BB^{N_j}$ for each
$j\le j_0$. Now, we must have
\begin{align*}
  \lambda(z, \bar z) ds^2_n-\sum_{j=1}^{j_0}\lambda_j(z, \bar z) F^*_j(ds^2_{N_j})=\sum_{j=j_0+1}^m\lambda_j(z, \bar z) F_j^*(ds^2_{N_j})
\end{align*}
in $(V_\gamma\cap \BB^n) \setminus S_{\gamma}$, which is a connected set by
Claim \ref{dimension}. Let $q\in (V_\gamma\cap \BB^n)\setminus
S_{\gamma}\rightarrow p\in \p \BB^n\cap V_\gamma$. Write
\begin{align*}\label{111}
(\lambda(z, \bar z)-\sum_{j=1}^{j_0} \lambda_j(z, \bar z))ds^2_n \bigg|_q +\sum_{j=1}^{j_0}
\lambda_j(z, \bar z)(ds^2_n-F^*_j(ds^2_{N_j}))\bigg |_q=\sum_{j=j_0+1}^m
\lambda_j(z, \bar z) F^*_j(ds^2_{N_j})\bigg |_q.
\end{align*}
 By Lemma 2.2, $X_j:=ds^2_n-F^*_j(ds^2_{N_j})$ is smooth up to $\p \BB^n$ for $j\le j_0$.
 We also see, by the choice of $j_0$ and Claim \ref{dimension}, that for a generic point $p$ in $\p \BB^n\cap V_{\gamma}$,
 $F^*_j(ds^2_{N_j})$ is real analytic in a small neighborhood of $p$ for each $j\ge j_0+1$.
 Thus by considering the normal component as before in the above equation, we see that
 $\lambda(z,\-{z})-\sum_{j=1}^{j_0}\lambda_j(z,\-{z})$
 vanishes to the order $\geq 2$ with respect to $1-|z|^2$ in an open set of the unit sphere. Since
 $\lambda(z,\-{z})-\sum_{j=1}^{j_0}\lambda_j(z,\-{z})$ is real
 analytic over $\CC^n$, we obtain
\begin{equation} \label{new-010}
 \lambda(z,\-{z})-\sum_{j=1}^{j_0} \lambda_j(z,\-{z})=(1-|z|^2)^2\psi(z,\-{z}).
\end{equation}
Here $\psi(z,\-{z})$ is a certain real analytic function over
$\CC^n$. Let

$$Y=(\lambda(z,\-{z})-\sum_{j=1}^{j_0}
 \lambda_j(z,\-{z}))ds^2_n.$$
Then $Y$ extends real analytically to $\CC^n$. Write
$X=\sum_{j=1}^{j_0} \lambda_j(z,\-{z}) X_j$. From what we argued above,
we easily see that there is  a certain small neighborhood $\mathcal
O$ of $q\in \p{\BB^n}$ in $\CC^n$ such that (1):  we can
holomorphically continue the initial germ of $F$ in $U$ through a
certain simple curve $\gamma$ with $\gamma(t)\in {\BB ^n}$ for $t\in
(0,1)$ to get a holomorpic map, still denoted by $F$, over $\mathcal
O$; (2): $\|F_j\|<1$ for
  $j>j_0$ and $\|F_j\|>1$ for $j\le j_0$ over ${\mathcal O\setminus \BB^n}$; and
(3):
\begin{equation}\label{new-100}
X=\sum_{j=1}^{j_0}
\lambda_j(z, \bar z)(ds^2_n-F^*_j(ds^2_{N_j})) \bigg( =\sum_{j=1}^{j_0}\lambda_j(z, \bar z)
X_j \bigg)=\sum_{j=j_0+1}^m \lambda_j(z, \bar z)F^*_j(ds^2_{N_j})-Y.
\end{equation}
We mention that we are able to make $\|F_j\|<1$ for any $z\in
{\mathcal O}$ and $j>j_0$ in the above due to the fact that
$(V_\gamma\cap {\BB ^n}) \setminus S_\gamma$, as defined before, is
connected.

\medskip

Now, let $\mathcal P$ be the union of the poles of $F_1, \ldots,
  F_{j_0}$. Fix a certain $p^*\in {\mathcal O}\cap \p {\BB^n}$ and
   let $\widetilde E=E\cup {\mathcal P}$. Then for any $\gamma: [0,1]\rightarrow \CC^n\setminus (\BB^n \cup \widetilde E)$ with $\gamma(0)=p^*$ and $\gamma(t)\not\in \p{\BB^n}$ for $t>0,\ \ F_j$ extends holomorphically to a small neighborhood $U_{\gamma}$ of
   $\gamma$ that contracts to $\gamma$.
   Still denote the holomorphic continuation of $F_j$ (from the initial germ of $F_j$ at $p^*\in {\mathcal O}$) over $U_{\gamma}$
    by $F_j$. If for some $t\in (0,1),
    \|F_j(\gamma(t))\|=1$, then we  similarly have

  \begin{claim}\label{dim10} Shrinking $U_\gamma$ if necessary, then
 $dim\big\{p\in U_{\gamma}: \|F_j(p)\|=1 \text{\ for some\ j} \big\} \le 2n-2$.
  \end{claim}

  {\noindent\bf Proof of the Claim \ref{dim10}:} Supposing otherwise we are going to deduce a contradiction. Define $S_\gamma$ in a similar way.
  Without loss of generality, we assume that $t_0\in (0,1)$ is the first point such that for a certain $j_{t_0}$,
   the local variety defined by $\|F_{j_{t_0}}(z)\|^2=1$ near $\gamma(t_0)$ has real dimension $2n-1$ at $\gamma(t_0)$.
   Then, as before, we have
\begin{align}\label{identity1}
X=\sum_{j=1}^{j_0}
\lambda_j(z, \bar z)(ds^2_n-F^*_j(ds^2_{N_j}))=\sum_{j=j_0+1}^{m} \lambda_j(z, \bar z)
F^*_j(ds^2_{N_j})-Y
\end{align}
in a connected component $W$ of $U_\gamma\setminus S_\gamma$ that
contains $\gamma (t)$ for $t<<1$ with $\gamma(t_0)\in \p W$.
Now, for any $q (\in W)\rightarrow p\in \p W$ near $p_0=\gamma(t_0)$
and $v\in {\CC^n}$ with $\|v\|=1$,
we have
the
following:
\begin{equation}\label{identity2}
\begin{split}
&\sum_{j=1}^{j_0} \lambda_j(q, \bar q) \bigg(ds^2_n(v, v)\bigg |_q - \frac{\|\sum_{\xi}\frac{\partial f_{j,l}}{\partial z_{\xi}}(q)v_{\xi}\|^2}{1-\|F_{j}(q)\|^2} - \frac{|\sum_{l,\xi}\overline {f_{j,l}(q)}\frac{\p f_{j,l}}{\p z_{\xi}}(q)v_{\xi}|^2}{(1-\|F_{j}(q)\|^2)^2}\bigg)\\
=&\sum_{j=j_0+1}^{m}  \lambda_j(q, \bar q)
\bigg(\frac{\|\sum_{\xi}\frac{\partial f_{j,l}}{\partial
z_{\xi}}(q)v_{\xi}\|^2}{1-\|F_{j}(q)\|^2}+\frac{|\sum_{l,\xi}\overline
{f_{j,l}(q)}\frac{\p f_{j,l}}{\p
z_{\xi}}(q)v_{\xi}|^2}{(1-\|F_{j}(q)\|^2)^2}\bigg)-Y(v, v)\bigg |_q.
\end{split}
\end{equation}
Now, if the local variety defined by $\|F_j(z)\|^2=1$ is not of real codimension one at $p_0$ for each $j\le j_0$, then the local variety
$S_{j'}$ defined by  $\|F_{j'}(z)\|^2=1$
has to be of real codimension one at $p_0$ for certain $j'> j_0$. Let
$J$ be the collection of all  such $j'$. Let $S^0$ be a small open piece of $\p W$ near $p_0$. Then for a generic $p\in
S^0$,
the left hand side of (\ref{identity2}) remains bounded as $q\ra
p\in S^0$. For a term on the right hand side with index $j\in J$,
if $S^0\cap S_{j}$ contains a germ of an irreducible component of $\p W$ of real codimension 1 containing
$p_0$, then it approaches to $+\infty$  for a generic $p$ unless
$F_{j}= constant$ as argued in the proof of Claim \ref{dimension}.
The other terms on the right hand side remain bounded as $q\ra p$
for a generic $p$. This is a contradiction to the assumption that
none of $F_j$'s for $j>j_0$ is constant.
 Hence, we can assume that the local variety defined by $\|F_j(z)\|^2=1$
near $p_0$ is of real codimension one for a certain $j\le j_0$.
Let $J$ be the set of indices such that for $j'\in J$, we have
$j'\le j_0$ and  $S_{j'}:=\{\|F_{j'}\|=1\}$ is the local real analytic
variety of real codimension one near $p_0$. For $j>j_0$, since
$\|F_{j}(z)\|<1$ for $z(\in U_\gamma)\approx p_0$ and since $t_0$ is
the first point such that $\|F_{j^*}\|=1$ defines a
variety of real codimension one for some $j^*$, we see that $\|F_{j}(z)\|<1$ for
$z(\in W)\approx p_0$. Define $S^0$ similarly, as an small open piece of $\p W$. Hence,
 as $q(\in W)\ra p\in
S^0$,
 the right hand side of (\ref{identity2}) is uniformly bounded from below.
  On the
other hand, on the left hand side of (\ref{identity2}), for any
$j'\in J$ with $S_{j'}\cap S^0$ containing an irreducible component of $\p W$ of real codimension 1 near $p_0$, if the numerator $|\sum_{l,\xi}\overline
{f_{j',l}(q)}\frac{\p f_{j',l}}{\p z_{\xi}}(q)v_{\xi}|^2$ of the
last term does not go to 0 for some vectors $v$, then the term with
index $j'$  on the left hand side would go to $-\infty$ for a
generic $p\in S^0$. If this happens to such $j'$, the left hand side
would approach to $-\infty$. Notice that all other terms on the
right hand side remain bounded from below as $q\ra p\in S^0$ for a generic $p$.
This is impossible. Therefore we must have for some $j'\in J$ that
$|\sum_{l}\overline {f_{j',l}(q)}\frac{\p f_{j',l}}{\p
z_{\xi}}(q)|^2=\frac{\p \sum_l |f_{j',l}|^2}{\p z_{\xi}}(q)=
\frac{\p \|F_{j'}\|^2}{\p z_{\xi}}(q)\rightarrow 0$ and thus
$\frac{\p \|F_{j'}\|^2}{\p z_{\xi}}(p)=0$ for all $\xi$ and $p\in
S_{j'}$. This immediately gives the equality $d(\|F_{j'}\|^2)=0$
along $S_{j'}$. Assume, without loss of generality, that $p_0$ is
also a smooth point of $S_{j'}$. If $S_{j'}$ has no complex
hypersurface passing through $p_0$, by a result of Trepreau [Tr],
the union of the image of local holomorphic disks attached to
$S_{j'}$ passing through $p_0$ fills in an open subset. Since
$F_{j'}$ is not constant, there is a small holomorphic disk smooth
up to the boundary $\phi(\tau): \BB^1\ra \CC^n$ such that $\phi(\p
\BB^1)\subset S_{j'}, \phi(1)=p_0$ and $F_{j'}$ is not constant
along $\phi$. Since $\p \BB^{N_{j'}}$ does not contain any
non-trival complex curves, $r=(\|F_{j'}\|^2-1)\circ \phi\not \equiv
0$.
 Applying the maximum principle and then the Hopf lemma to the
subharmonic function $r=(\|F_{j'}\|^2-1)\circ \phi$, we see that the
outward normal derivative of $ r$  at $\tau=1$ is positive.
This contradicts to the fact that $d(\|F_{j'}\|^2)=0$ along
$S_{j'}$. We can argue in the same way for points $p\in S_{j'}$ near
$p_0$ to conclude that for any $p\in S_{j'}$ near $p_0$, there is a
complex hypersurface contained in $S_{j'}$ passing through  $p$.
Namely, $S_{j'}$ is  Levi flat, foliated by a family of smooth
complex hypersurfaces denoted by  ${Y_\eta}$ with real parameter
$\eta$ near $p_0$. Let $Z$ be a holomorphic vector field along
$Y_\eta$. We then easily see that
$0=\overline{Z}Z(\|F_{j'}\|^2-1)=\sum_{k=1}^{N_{j'}}|Z(f_{j',k})|^2$.
Thus, we see that $F_{j'}$ is constant along each $Y_\eta$. Hence,
$F_{j'}$ cannot be a local embedding at each point of $S_{j'}$. However, on
the other hand, notice that $F_{j'}$ is a proper holomorphic map
from $\BB^n$ into ${\BB}^{N'}$, then $F_{j'}$ is a local embedding
near $\p\BB^n$. This implies that the set of points where $F_{j'}$ is not a
local embedding can be at most of complex codimension one (and thus
real codimension two). This is a contradiction. This proves Claim
\ref{dim10}.

\medskip
Hence, we see that ${\mathcal E}=\{p\in \CC^n\setminus
(\overline{\BB^n}\cup \widetilde E): $ some branch, obtained by
the holomorphic continuation through curves described before, of $\
F_j \ \text{for some}\  j\ \text{maps} \ p \ \text{to} \ \p
\BB^{N_j}\}$
 is a real analytic variety of real dimension at most $2n-2$. Now,
  for any $p\in \CC^n\setminus (\overline{\BB^n}  \cup {\widetilde E})$, any curve
  $\gamma: [0,1]\rightarrow \CC^n\setminus (\BB^n \cup
 \widetilde E)$
  with $\gamma(0)=p^*\in {\mathcal O}\cap \p{\BB^n}, \gamma(t)\not\in \p{\BB^n}$ for $t>0$
  and $\gamma(1)=p$,
  we can homotopically change $\gamma $ in $\CC^n\setminus (\overline{\BB^n}\cup \widetilde
  E)$ (but without changing the terminal point)
  such that $\gamma(t)\notin {\mathcal E}$ for $t\in (0,1)$. Now,
  the holomorphic continuation of the initial germ of $F_j$ from
  $p^*$
  never cuts $\p \BB^{N_j}$ along $\gamma(t)$ ($0< t< 1$). We thus see that $\|F_j(p)\|\le 1$ for $j>j_0$.

Let $\{(f_{j,l})_{k;p}\}_{k=1}^{n_{jl}}$ be  all possible (distinct)
germs of holomorphic functions that we can get at $p$ by the
holomorphic continuation, along curves described above in
$\CC^n \setminus (\BB^n \cup {\tilde E})$,  of
$f_{j,l}$. Let $\sigma_{jl,\tau}$ be the fundamental symmetric
function of $\{(f_{j,l})_{k;p}\}_{k=1}^{n_{jl}}$ of degree $\tau$.
Then $\sigma_{jl, \tau}$ well defines a holomorphic function over
$\CC^n\setminus \overline{\BB^n} $. $|\sigma_{jl,\tau}|$ is bounded in
$\CC^n\setminus (\overline{\BB^n}\cup {\tilde E})$. By the
Riemann removable singularity theorem, $\sigma_{jl, \tau}$ is
holomorphic over $\CC^n\setminus\overline{\BB^n}$. By the Hartogs
extension theorem, $\sigma_{jl, \tau}$ extends to a  bounded
holomorphic function over $\CC^n$.
 Hence, by the Liouville theorem, $\sigma_{jl, \tau}\equiv const$.
This forces $(f_{j,l})_{k}$ and thus $F_j$ for $j>j_0$  to be
constant. We obtain a contradiction. This proves that each $F_j$
extends to a proper rational map from ${\BB^n}$ into $\BB^{N_j}$.
Together with Proposition \ref{mainprop}, we complete the proof of
the main Theorem. $\endpf$


\begin{remark}The regularity of
$\lambda_j,\lambda$ can be reduced to be only real analytic in the
complement of a certain real codimension two subset. Also, we need
only to assume that they are positive outside a real analytic
variety of real codimension two. This is obvious from our proof of
Theorem 1.1.
\end{remark}

\begin{remark}
Assume that $\lambda, \lambda_j$ are smooth, positive, Nash algebraic (or more generally, real analytic) functions on $\mathbb{B}^n, \mathbb{B}^{N_j}$ respectively for all $j$ and $F=(F_1, \cdots, F_m): U \subset \mathbb{B}^n \rightarrow \mathbb{B}^{N_1} \times \cdots \times \mathbb{B}^{N_m}$ is a holomorphic embedding such that
$$\lambda ds^2_n = \sum_{j=1}^m F^*_j(\lambda_j ds^2_{N_j}).$$
It would be very interesting to prove the total geodesy for each non-constant component $F_j$. However, different from the situation in Theorem \ref{main1}, one cannot prove the algebraic extension of $F$ using the technique in the appendix since we do not know yet how to construct a target real algebraic hypersurface associated to $F$. Once the algebraic extension of $F$ is obtained, the total geodesy should follow from our argument without much modification. For the related algebraic extension problem, see \cite{HY}.
\end{remark}


\section{Appendix: algebraic extension}
In this appendix, we prove the algebraicity of the local holomorphic map $F$ in Theorem \ref{main1}.
As in the theorem, we let $U\subset \BB^n$ be a
connected  open subset. Let
\begin{align*}
&F=(F_1,\ldots, F_m): U
\rightarrow \BB^{N_1}\times
\cdots \times \BB^{N_m}
\end{align*}
be a holomorphic isometric embedding up to conformal factors
$\{\lambda(z,\-{z});
\lambda_1(z,\-{z}),\cdots,\lambda_m(z,\-{z})\}$ in the sense that
 $$\lambda(z,\-{z}) ds^2_n=\sum_{j=1}^{m}\lambda_j(z,\overline{z})
F_j^*(ds^2_{N_j}).$$ Here $\lambda_j(z,
\bar z)>0,  \lambda(z, \bar z)>0$ are smooth Nash algebraic
functions in $\CC^n$, and $ds_{n}^2$ and $ds_{N_j}^2$ are the
Bergman metrics of $\BB^n$ and $\BB^{N_j}$, respectively. We further
assume without loss of generality that none of the $F_j$'s is a
constant map. Our proof uses  exactly the same method employed  in
the paper of Mok-Ng \cite{MN}, following a suggestion of Yum-Tong
Siu. Namely, we use the Grauert tube technique to reduce the problem
to the algebraicity problem for CR mappings.
However, different from the consideration in \cite{MN}, the Grauert
tube constructed by using the unit sphere bundle over
$\BB^{N_1}\times \cdots \times \BB^{N_m}$ with respect to the metric
$\oplus_{j=1}^m ds_{N_j}^2$, up to  conformal factors, may have
complicated geometry and may not even be pseudoconvex anymore in
general.
 To overcome the difficulty,  we bend the target hypersurface to make it sufficiently positively
curved along the tangential direction of the source domain. 

\medskip

Let $K>0$ be a large constant to be determined. Consider $S_1 \subset T U$ and $S_2 \subset U \times T
\BB^{N_1} \times \cdots \times T \BB^{N_m}$ as follows:

\begin{equation}
S_1 := \left\{(t, \zeta)  \in T U : (1+K |t |^2)\lambda(t, \bar
t)ds^2_{n}(t)(\zeta, \zeta) =1 \right\},
\end{equation}

\begin{equation}
\begin{split}
S_2 &:=\{ (t, z_1, \xi_1, \cdots, z_m, \xi_m) \in U \times T \BB^{N_1} \times \cdots \times T \BB^{N_m}:\\
&\ \ \ \ (1+K |t |^2)[\lambda_1(t, \bar t)ds^2_{N_1}(z_1)(\xi_1, \xi_1) + \cdots +
\lambda_m(t, \bar t)ds^2_{N_m}(z_m)(\xi_m, \xi_m)]=1\}.
\end{split}
\end{equation}
The defining functions $\rho_1, \rho_2$ of $S_1, S_2$ are,
respectively, as follows:
$$\rho_1=(1+K |t |^2)\lambda(t, \bar t)ds^2_n(t)(\zeta, \zeta)-1,$$
$$\rho_2=(1+K |t |^2)[\lambda_1(t, \bar t)ds^2_{N_1}(z_1)(\xi_1, \xi_1) + \cdots + \lambda_m(t, \bar t)ds^2_{N_m}(z_m)(\xi_m, \xi_m) ] - 1.$$
Then one can easily check that the map $(id, F_1, dF_1, \cdots, F_m, dF_m)$
maps $S_1$ to $S_2$ according to the metric equation

$$\lambda(t,\-{t}) ds^2_n=\lambda_1(t,\-{t})
F_1^*(ds^2_{N_1}) + \cdots+ \lambda_m(t,\-{t}) F_m^*(ds^2_{N_m}).$$

\medskip

\begin{lemma}\label{pseudoconvex}
$S_1, S_2$ are both real algebraic hypersurfaces. Moreover for $K$
sufficiently large, $S_1$ is smoothly strongly pseudoconvex. For any
$\xi_1\not=0, \cdots, \xi_m\not=0$, $(0, 0 , \xi_1, \cdots , 0, \xi_m) \in S_2$ is a smooth
strongly pseudoconvex point when $K$ is sufficiently large, where
$K$ depends on the choice of $\xi_1, \cdots, \xi_m$. 
\end{lemma}

{\noindent \bf Proof of Lemma \ref{pseudoconvex}:}
It is  immediate from the defining functions that $S_1, S_2$ are
smooth real algebraic hypersurfaces. We show the strong
pseudoconvexity of $S_2$ at $(0, 0, \xi_1, \cdots, 0, \xi_m)$ as follows: (The
strong pseudoconvexity of $S_1$ follows from the same computation.)

By applying $\p\bar\p$ to $\rho_2$ at $(0, 0, \xi_1, \cdots,  0, \xi_m)$, we
have the following Hessian matrix
\begin{equation}
\begin{bmatrix}
A & 0 & D_1 & \cdots & 0 & D_m\\
0 & B_1 & 0 & \cdots & 0 & 0\\
\bar D_1 & 0 & C_1 &\cdots & 0 & 0 \\
\cdots & \cdots & \cdots & \cdots & \cdots & \cdots \\
0 & 0 & 0 & \cdots & B_m & 0 \\
\bar D_m & 0 & 0 & \cdots & 0 & C_m \\
\end{bmatrix}
\end{equation}
where $A, B_j, C_j, D_j, j=1,2, \cdots, m$ are function-valued matrices with the following (in)equalities:
\begin{equation}
\begin{split}
A:&=\bigg(\p_{t_i}\p_{t_{\bar j}} \rho_2\bigg)= \bigg(K \sum_{\nu=1}^m \lambda_\nu(0)|\xi_\nu|^2 \delta_{ij}+ \sum_{\nu=1}^m \p_{t_i}\p_{t_{\bar j}}\lambda_\nu(0)|\xi_\nu|^2 \bigg) \\
&\geq \delta K (|\xi_1|^2+\cdots + |\xi_m|^2)I_n,
\end{split}
\end{equation}
\begin{equation}
\begin{split}
B_j:= \bigg(\p_{z_{jk}}\p_{\bar z_{
{jl}}} \rho_2\bigg)= \bigg(-\lambda_j(0) \sum_{\nu,\mu=1}^{N_j} R_{z_{jk} \bar z_{jl} \mu
\bar{\nu}}(0) \xi_{j\mu} \bar\xi_{j\nu}\bigg) \geq \delta |\xi_j|^2 I_{N_j},
\end{split}
\end{equation}
\begin{equation}
\begin{split}
C_j:=\bigg(\p_{\xi_{jk}}\p_{\bar \xi_{j l}}\rho_2\bigg)=\bigg(\lambda_j(0)\delta_{k \bar l}\bigg) \geq \delta
I_{N_j}, 
\end{split}
\end{equation}
\begin{equation}
\begin{split}
D_j:=\bigg(\p_{t_i}\p_{\bar \xi_{j l}}\rho_2\bigg)= \bigg(\p_{t_i}\lambda_j(0) \xi_{jl}\bigg)_{i\le n}^{l\le N_j},
\end{split}
\end{equation} at $(0, 0, \xi_1, \cdots, 0, \xi_m) $ for some $\delta>0$.

Let $(e, r_1, s_1, \cdots, r_m, s_m)
\not= 0$, where $e=(e_1, \cdots, e_n), r_j=(r_{j1}, \cdots r_{j N_j}), s_j=(s_{j1}, \cdots s_{j N_j})$ for all $j$. It holds that
\begin{equation}
\begin{split}
~~& \begin{bmatrix}e & r_1 & s_1 &\cdots & r_m & s_m
\end{bmatrix}
\begin{bmatrix}
A & 0 & D_1 & \cdots & 0 & D_m\\
0 & B_1 & 0 & \cdots & 0 & 0\\
\bar D_1 & 0 & C_1 &\cdots & 0 & 0 \\
\cdots & \cdots & \cdots & \cdots & \cdots & \cdots \\
0 & 0 & 0 & \cdots & B_m & 0 \\
\bar D_m & 0 & 0 & \cdots & 0 & C_m \\
\end{bmatrix}
\begin{bmatrix}\bar e^t \\ \bar r^t_1 \\ \bar s^t_1 \\  \cdots \\ \bar r_m^t \\ \bar s_m^t
\end{bmatrix} \\
&\geq \delta K |e|^2 \sum_{j=1}^m |\xi_j|^2 + \delta \sum_{j=1}^m |\xi_j|^2|r_j|^2 +\delta \sum_{j=1}^m |s_j|^2 
-2 \sum_{j=1}^m \left|\sum_{i\le n,\ l\le N_j} e_i \p_{t_i}\lambda_j(0)\xi_{jl} \bar s_{jl}
\right| 
\\
&\geq \sum_{j=1}^m \left[ (\delta K-M) |\xi_j|^2 |e|^2 +  \delta |\xi_j|^2|r_j|^2 + (\delta - \epsilon)|s_j|^2 \right]\\
&> 0.
\end{split}
\end{equation}
Here the second inequality holds since

$$\left|\sum_{i,l}
e_i \p_{t_i}\lambda_j(0)\xi_{jl} \bar s_{jl} \right| \leq M_1|e||\xi_j
\cdot \bar s_j| \leq M_1 |e||\xi_j||s_j| \leq \frac{M}{2}|e|^2|\xi_j|^2 +
\frac{\epsilon}{2}|s_j|^2$$ by the standard Cauchy-Schwarz inequality and
$M=\frac{M_1^2}{\epsilon}$. The last strict inequality holds as $\xi_j
\not =0$ for all $j$ 
by letting $\epsilon < \delta$ and
raising $K$ sufficiently large. $\endpf$


\begin{theorem}\label{Nash}
Under the assumption of Theorem 1.1, $F$ is Nash algebraic.
\end{theorem}

{\noindent \bf Proof of Theorem \ref{Nash}:} Without loss of generality, one can assume that $F(0)=0$ by
composing elements from $Aut(\BB^n)$ and $Aut(\BB^{N_1})\times \cdots \times
Aut(\BB^{N_m})$. Furthermore, since $F_1, \cdots, F_m$ are not constant
maps, we can assume  that $d F_1 |_0\not\equiv  0, \cdots, d F_m
|_0\not\equiv 0$. Therefore, there exists $0 \not= \zeta \in T_0
\BB^n$, such that $d F_j(\zeta) \not =0$ for all $j$. 
After scaling, we assume that $(0,\zeta)\in S_1$. Notice that both
the fiber of $S_1$ over $0\in U$ and the fiber of $S_2$ over
$(0,0, \cdots, 0)\in U\times {\BB^{N_1}}\times \cdots \times {\BB^{N_m}}$ are independent of
the choice of $K$.
 Now the theorem follows by applying the algebracity theorem
of Huang [Hu1] and Lemma \ref{pseudoconvex} to the map
$(id, F_1, dF_1, \cdots, F_m, d F_m)$ from $S_1$ into $S_2$. $\endpf$


\end{document}